\documentclass[a4paper,11pt]{article}
\usepackage[russian,english]{babel}
\usepackage{amsfonts,amsmath,amssymb,theorem}
\usepackage{hyperref} % hyperlinks
\usepackage{colortbl}
\usepackage{ifthen}
\usepackage{hhline}

\makeatletter
\def\@seccntformat#1{\csname the#1\endcsname.\ } % точка после номера раздела
\def\@biblabel#1{#1.} % формат номеров в списке литературы
\makeatother

\date{}

    \newif\ifNoRemark
    \def\addtheorem#1#2#3#4{
    \ifthenelse{\expandafter\isundefined\csname the#2\endcsname}{\newcounter{#2}}{}
    \newenvironment{#1}[1][\global\NoRemarktrue]% No Remark by default
     {\par\addvspace{2mm plus 0.5mm minus 0.2mm}\noindent
       \refstepcounter{#2}{\bf #3~\csname the#2\endcsname
      \vphantom{##1}\ifNoRemark.\ \else\ (##1).\fi}\begingroup #4}%
      {\endgroup\par\addvspace{1mm plus 0.5mm minus 0.2mm}\global\NoRemarkfalse}
    \expandafter\newcommand\csname b#1\endcsname{\begin{#1}}
    \expandafter\newcommand\csname e#1\endcsname{\end{#1}}
    }

\addtheorem{theorem}{thrm}{Theorem}{\sl}
\addtheorem{lemma}{lmm}{Lemma}{\sl}
\addtheorem{pro}{prpstn}{Proposition}{\sl}
\addtheorem{corol}{corol}{Corollary}{\sl}
\begin{document}

\title{On
 perfect $2$-colorings of the $q$-ary $n$-cube  \thanks{The work is supported by RFBR (grants
10-01-00424, 10-01-00616) and Federal Target Grant "Scientific and
educational personnel of innovation Russia" for 2009-2013
(government contract No. 02.740.11.0429)} }
\author{Vladimir N.~Potapov}

\maketitle

\begin{center}
\textit{Sobolev Institute of Mathematics,\\
Novosibirsk State University, Novosibirsk }
\end{center}

\begin{abstract}

A coloring of the $q$-ary $n$-dimensional cube (hypercube) is called
perfect if, for every $n$-tuple $x$,  the collection of the colors
of the neighbors of $x$ depends only on the color of $x$.  A
Boolean-valued function is called correlation-immune of degree $n-m$
if it takes the value $1$ the same number of  times for each
$m$-dimensional face of the hypercube. Let $f=\chi^S$ be a
characteristic function of some subset $S$ of  hypercube. In the
present paper it is proven  the inequality $\rho(S)q({\rm
cor}(f)+1)\leq \alpha(S)$,
   where ${\rm cor}(f)$
  is the maximum degree of the correlation immunity of $f$, $\alpha(S)$
  is the average number of neighbors in the set $S$
for $n$-tuples in the complement of a set $S$, and
  $\rho(S)=|S|/q^n$ is the density of the set $S$.
  Moreover, the function $f$ is a perfect coloring if and only
  if we obtain an equality in the above formula.
Also we find new lower bound for the cardinality of components of
    perfect coloring and 1-perfect code in the case $q>2$.

{\bf Keywords}: hypercube, perfect coloring, perfect code, MDS code,
bitrade, equitable partition,  orthogonal array.

\end{abstract}

\section{Introduction}

Let $Z_q$ be the set of entries $\{0,\dots,q-1\}$. The set $Z_q^n$
of $n$-tuples of entries is called $q$-ary $n$-dimensional cube
(hypercube).  The {\it Hamming distance} $d({x},{y})$ between two
$n$-tuples ${x}, {y}\in Z_q^n$ is the number of positions at which
they differ. Define the  number $\alpha(S)$ to be the average number
of neighbors in the set $S\subseteq Z_q^n$ for $n$-tuples in the
complement of a set $S$, i.\,e., $\alpha(S)=
\frac{1}{q^n-|S|}\sum\limits_{x\not \in S}|\{y\in S \ |\
d(x,y)=1\}|$.

A mapping $Col:Z_q^n\rightarrow\{0,\dots,k\}$ is called a {\it
perfect
 coloring} with matrix of parameters $A=\{a_{ij}\}$  if, for all $i$,
  $j$, for every  $n$-tuple of color $i$, the number of  its neighbors of color $j$ is
  equal to $a_{ij}$. Other terms used for this notion in the literature are
"equitable partition", "partition design" and "distributive
coloring". In what follows we will only consider colorings in two
colors (2-coloring). Moreover, for convenience we will assume that
the set of colors is $\{0,1\}$. In this case the  Boolean-valued
function $Col$ is a characteristic  function of the set of 1-colored
$n$-tuples.

   A {\it 1-perfect code} (one-error-correcting) $C\subset Z_q^n$ can be defined
    as the set of units of a perfect coloring with
   matrix of parameters $  A= \left(
 \begin{array}{ccc}
n(q-1)-1 &  1 \\
  n(q-1)  & 0 \\
  \end{array} \right).
$ If $q$ is the power of a prime number then the coloring with such
parameters exists only if $n=\frac{q^m-1}{q-1}$ ($m$ is integer).
For $q=2$ a list of accessible parameters and corresponding
constructions of perfect 2-colorings can be found in \cite{Fl1} and
\cite{Fl3}.

  In \cite{Fl2} it is established  that for each unbalanced Boolean
  function $f=\chi^S$ ($S\subset Z_2^n$)
  the inequality
  ${\rm cor}(f)\leq \frac{2n}{3} -1$  holds. Moreover, in the case of the equality
    ${\rm cor}(f)= \frac{2n}{3} -1$, the function $f$  is a perfect 2-coloring.
Similarly, if for any set $S\subset  Z_2^n$ the Friedman (see
\cite{Fr}) inequality  $\rho(S)\geq 1-\frac{n}{2({\rm cor}(f)+1)}$
becomes an equality then the function $\chi^S$ is a perfect
$2$-coloring (see \cite{Pot}). Consequently, in the extremal cases,
the regular distribution on balls follows from the uniform
distribution on faces. The main result of present paper is following
theorem: \btheorem\label{thki1}

 {\rm (a)} For each Boolean-valued function
   $f=\chi^S$, where $S\subset Z_q^n$, the inequality\\
   $\rho(S)q({\rm
cor}(f)+1)\leq \alpha(S)$   holds.

 {\rm (b)} A Boolean-valued function $f=\chi^S$ is  a perfect $2$-coloring
 if and only if\\
$\rho(S)q({\rm cor}(f)+1)=\alpha(S)$. \etheorem

\section{Criterion for perfect 2-coloring}

In the proof of the theorem we employ the idea from the papers
\cite{Br}.

Now we consider $Z_q$ as the cyclic group on the set of entries
$\{0,\dots,q-1\}$. We may impose the structure of the group
$Z_q\times\dots \times Z_q$ on the hypercube. Consider the vector
space $\mathbb{V}$ of complex-valued function on $Z_q^n$ with scalar
product $(f,g)=\frac{1}{q^n}\sum\limits_{x\in
Z_q^n}f(x)\overline{g(x)}$. For every $z\in Z_q^n$ define a {\it
character} $\phi_z(x)=\xi^{\langle x, z\rangle}$, where $\xi=
e^{2\pi i/q}$ is a primitive complex $q$th root of unity and
${\langle x, z\rangle}= x_1z_1+\dots+x_nz_n$. Here all arithmetic
operations are performed on   complex numbers. As is generally known
the characters of the group $Z_q\times\dots \times Z_q$ form an
orthonormal basis of $\mathbb{V}$. It is sufficient to verify that
$\xi^k\overline{\xi^k} =1$ and $\sum\limits_{j=0}^{q-1}\xi^{kj}=0$
as $k\neq 0\mod q$.

Let $M$ be the adjacency matrix by the hypercube $Z_q^n$. This means
that $Mf(x)=\sum\limits_{y,d(x,y)=1}f(y)$. It is well known that the
characters are eigenvectors of $M$. Indeed we have
$$M\phi_z(x)= \sum\limits_{y,d(x,y)=1}\xi^{\langle y-x, z\rangle+\langle x, z\rangle}=
\xi^{\langle x, z\rangle}\sum\limits_{j=1}^n \sum\limits_{k\neq
0}\xi^{kz_j} = ((n-wt(z))(q-1) - wt(z))\phi_z(x), $$ where $wt(z)$
is the number of nonzero coordinates of $z$.

Consider a perfect coloring $f\in \mathbb{V}$, $f(Z_q^n)=\{0,1\}$
with matrix of parameters
\begin{equation}\label{echar1} A=\left(
 \begin{array}{ccc}
n(q-1)-b & b  \\
  c  & n(q-1)-c  \\
  \end{array} \right).
\end{equation}
   The vector $(-b,c)$  is an eigenvector of
  $A$ with the eigenvalue $n(q-1) -c -b$. The definition of a perfect
  2-coloring implies that the function $(b+c)f-b$ is the eigenvector of matrix
  $M$. Moreover the converse is true: every two-valued eigenvector of matrix
  $M$  generates a perfect coloring.

\bpro\label{stki2} (see \cite{Fl1})

{\rm (a)} Let $f$ be a perfect 2-coloring with matrix of parameters
$A$ (\ref{echar1}). Then $s=\frac{c+b}{q}$ is integer and
$(f,\phi_z)=0$ for every  $n$-tuples $z\in Z_q^n$ such that
$wt(z)\neq 0,s$.

{\rm (b)} Let $f:Z_q^n\rightarrow \{0,1\}$ be a Boolean-valued
function. If $(f,\phi_z)=0$ for every $n$-tuples $z\in
\{0,\dots,q-1\}^n$ such that $wt(z)\neq 0,s$ then $f$ is a perfect
2-coloring. \epro

Refer as a {\it correlation-immune} function of order $n-m$ to a
function $f\in \mathbb{V}$ that every value is uniformly distributed
on all $m$-dimensional faces. For any function $f\in \mathbb{V}$ we
denote the maximum of order of its correlation-immunity by ${\rm
cor}(f)$. Consider a nonempty set of $n$-tuples $O(a)=
f^{-1}(a)\subset Z^n_q$ where $a\in \mathbb{C}$. An array consisted
of $n$-tuples $x\in O(a)$ is called {\it orthogonal array} with
parameters $OA_\lambda({\rm cor}(f),n,q)$, where
$\lambda=|O(a)|/q^{n-{\rm cor}(f)}$.

\bpro\label{prochar2} (see \cite{Br})

{\rm (a)} If $f\in \mathbb{V}$ is a correlation-immune function of
order $m$ then $(f,\phi_z)=0$ for every $n$-tuples $z\in Z_q^n$ such
that $0<wt(z)\leq m$.

{\rm (b)} A Boolean-valued function $f\in \mathbb{V}$ is
correlation-immune  of order $m$ if $(f,\phi_z)=0$ for every
$n$-tuples $z\in Z_q^n$ such that $0<wt(z)\leq m$. \epro

\bcorol\label{corchar}
 Let $f$ be a perfect 2-coloring with matrix of parameters
(\ref{echar1}). Then ${\rm cor}(f)=\frac{c+b}{q}-1$. \ecorol

For 1-perfect codes last statement was proven otherwise in
\cite{Del}.

{\bf Proof of the theorem.} We have the following equalities by the
definitions and general properties of orthonormal basis.

\begin{equation}\label{echar3}
\sum\limits_z|(f,\phi_z)|^2=\frac{1}{q^n}\sum\limits_{x\in
Z^n_q}|f(x)|^2= \rho(S).
\end{equation}

\begin{equation}\label{echar31}
(f,\phi_{\overline{0}})=\frac{1}{q^n}\sum\limits_{x\in Z^n_q}f(x)=
\rho(S).
\end{equation}

\begin{equation}\label{echar4}
(Mf,f)=\frac{1}{q^n}\sum\limits_{x\in
Z^n_q}\sum\limits_{y,d(x,y)=1}f(x)\overline{f(y)}= {\rm nei}
(S)\rho(S),
\end{equation}

where ${\rm nei} (S)=\frac{1}{|S|}\sum\limits_{x \in S}|\{y\in S \
|\ d(x,y)=1\}|$.

\begin{equation}\label{echar5}
(Mf,f)= \sum\limits_{z\in Z^n_q}(n(q-1) - wt(z)q)|(f,\phi_z)|^2.
\end{equation}

  From (\ref{echar3}--\ref{echar5}) and  Proposition \ref{prochar2} we obtain the
equality
$${\rm nei}(S)\rho(S)= \rho(S)^2n(q-1) +
\sum\limits_{z,wt(z)\geq {\rm cor}(f)+1}(n(q-1) -
wt(z)q)|(f,\phi_z)|^2.$$
 Since  $\sum\limits_{z,wt(z)\geq {\rm cor}(f)+1}|(f,\phi_z)|^2= \rho(S) -
 \rho(S)^2$,
 we have
$${\rm nei}(S)\rho(S)\leq \rho(S)^2n(q-1) +
(n(q-1)-({\rm cor}(f)+1)q)(\rho(S) -
 \rho(S)^2)\ \rm{and}$$
\begin{equation}\label{echar10}
 ({\rm
cor}(f)+1)q(1 -
 \rho(S))\leq n(q-1)-{\rm nei}(S).
\end{equation}
Substitute the set $Z^n_q\setminus S$  instead of the set $S$ into
the inequality (\ref{echar10}). Since ${\rm cor}(\chi^{S})={\rm
cor}(\chi^{Z^n_q\setminus S})$, $1 -
 \rho(Z^n_q\setminus S)=\rho(S)$ and
  $n(q-1)-{\rm nei}(Z^n_q\setminus S)=\alpha(S)$ we obtain
the item (a) of the Theorem.

 Moreover, the equality
\begin{equation}\label{ki17}
({\rm cor}(f)+1)q(1 -
 \rho(S))= n(q-1)-{\rm nei}(S)
\end{equation}
 holds if and only if $(f,\phi_z)=0$ for every $n$-tuple $z$ such that
$wt(z)\geq {\rm cor}(f) +2$. Then from Proposition \ref{stki2} (b)
we conclude that $f$ is a perfect $2$-coloring.

Each perfect $2$-coloring satisfies (\ref{ki17}), which is a
consequence of Proposition \ref{stki2} (a) and  Corollary
\ref{corchar}. As mentioned above  the equality (\ref{ki17}) is
equivalent to the equality in the item (b) of the Theorem. $\Box$

Since ${\rm nei}(S)\neq 0$, the inequality (\ref{echar10}) implies
the Bierbrauer -- Friedman inequality (see \cite{Fr}, \cite{Br})
$$\rho(S)\geq 1-\frac{n(q-1)}{q({\rm cor}(f) +1)}.$$

For 1-perfect binary codes,  a similar  theorem was previously
proven in \cite{OP}. Namely, if  ${\rm cor}(S)={\rm cor}(H)$ and
  $\rho(S)=\rho(H)$, where
  $S,H\subset Z_2^n$ and  $H$  is a 1-perfect code, then $S$ is also a 1-perfect code.

\section{Components of perfect 2-coloring}

Refer as a {\it bitrade}  of order $n-m$ to a subset  $B\subseteq
Z^n_q$ that the cardinality of intersections $S$ and each
$m$-dimensional face are even.

\bpro\label{prochar3} Let $S\subseteq Z^n_2$ be a nonempty {\it
bitrade} of order $m$. Then $|S|\geq 2^{m+1}$.  \epro

Proposition \ref{prochar3} formulated  in other term was proven in
\cite{MW}.

\bpro\label{prochar4} Let $S\subseteq Z^n_q$ ($q>2$) be a nonempty
{\it bitrade} of order $m$. Then $|S|\geq 2^{m+1}$.  \epro {\bf
Proof.} Suppose that this statement is true for $n=k$. We will prove
it for $n=k+1$. Let there exist three parallel $k$-dimensional faces
$F_1$, $F_2$, $F_3$ such that intersections $F_i\cap S$ are
nonempty. By induction hypothesis $|F_i\cap S|\geq 2^{m-1}$ for
$i=1,2,3$; consequently, $|S|\geq 3\cdot 2^{m-1}$. In the other case
$|S|\geq 2^m$ by Proposition \ref{prochar3}. $\Box$

Let  characteristic functions $f=\chi^{S_1}$ and $g=\chi^{S_2}$ be
perfect 2-colorings (correlation-immune) with an equal matrix of the
parameters (${\rm cor}(f)={\rm cor}(g)$). A set $S_1\bigtriangleup
S_2$ is called {\it mobile} and  sets $S_1\setminus S_2$ and
$S_2\setminus S_1$ are called {\it components} of a perfect
2-colorings (correlation-immune functions) $\chi^{S_1}$ and
$\chi^{S_2}$ respectively. It is clear, that a mobile set of
correlation-immune function of order $m$ is a bitrade of order $m$.

\bcorol\label{corchar2}

(a)  Let $f$ be a perfect 2-coloring with matrix of parameters
(\ref{echar1}). If $S\subset Z^n_q$ is a component of $f$ then
$|S|\geq 2^{\frac{c+b}{q}-1}$.

b)  Let $C\subset Z^n_p$ be a 1-perfect code. If $S\subset Z^n_q$ is
a component of $f$ then $|S|\geq 2^{\frac{n(q-1)+1}{q}-1}$. \ecorol

If $q=2$ then the lower bound $|S|\geq 2^{\frac{n+1}{2}-1}$ for the
cardinality of components of 1-perfect codes is achievable (see, for
example, \cite{Pot}). In the case $q>2$  an upper bound  for the
cardinality of components of 1-perfect codes is obtained
constructively (see \cite{PV}, \cite{Los}). If $q=p^r$ and $p$ is a
prime number then $|S|\geq p^{\frac{q^{m-1}-1}{q-1}(r(q-2)+1)}$
where $n=\frac{q^m-1}{q-1}$.

A set $S\subset Z^n_p$ is called {MDS code} with distance $2$ if
 intersection $S$ and each 1-dimensional face contains precisely  one
$n$-tuple. Obviously a characteristic function of  MDS code is a
perfect 2-coloring with matrix of parameters $\left(
 \begin{array}{ccc}
n(q-2) &  n \\
 n(q-1)   &  0 \\
  \end{array} \right)$. If $q\geq 4$ then the lower bound $|S|\geq 2^{n-1}$ for the
cardinality of  the components of MDS codes is achievable (see
\cite{PKS}).

\begin{center}

\end{center}

\end{document}